\documentclass{crm-eca}

\usepackage{amssymb}

\usepackage{Irtatca}  

\newtheorem{theorem}{Theorem}

\newtheorem{proposition}[theorem]{Proposition}

\theoremstyle{definition}

\newtheorem{example}[theorem]{Example}

\newtheorem{question}[theorem]{Question}
\newtheorem{open}[theorem]{Open Problem}
\newtheorem{partial answer}[theorem]{Partial Answer}

\newtheorem*{aim}{Aim}
\newtheorem*{notation}{Notation}
\newtheorem{DefProp}[theorem]{Definition/Proposition}
\newtheorem*{proposition*}{Proposition}

\theoremstyle{remark}
\newtheorem{remark}[theorem]{Remark}

\newcommand{\opname}[1]{\operatorname{\mathsf{#1}}}

\newcommand{\ca}{{\mathcal A}}
\newcommand{\cb}{{\mathcal B}}
\newcommand{\cc}{{\mathcal C}}
\newcommand{\cd}{{\mathcal D}}
\newcommand{\ce}{{\mathcal E}}
\newcommand{\cf}{{\mathcal F}}
\newcommand{\cp}{{\mathcal P}}
\newcommand{\cm}{{\mathcal M}}
\newcommand{\ct}{{\mathcal T}}

\renewcommand{\mod}{\opname{mod}\nolimits}
\newcommand{\ul}[1]{\underline{#1}}
\newcommand{\proj}{\opname{proj}\nolimits}
\newcommand{\add}{\opname{add}\nolimits}

\renewcommand{\ker}{\opname{ker}\nolimits}
\DeclareMathOperator{\End}{\mathsf{End}}
\DeclareMathOperator{\Hom}{\mathsf{Hom}}
\DeclareMathOperator{\Ext}{\mathsf{Ext}}
\DeclareMathOperator{\krdim}{\mathsf{kr.dim}}

\DeclareMathOperator{\injdim}{\mathsf{inj.dim}}
\DeclareMathOperator{\gldim}{\mathsf{gl.dim}}
\DeclareMathOperator{\virdim}{\mathsf{vir.dim}}
\def\GP{\mathop{\rm GP}\nolimits}

\begin{document}

\papertitle{A remark on Leclerc's Frobenius categories }


\paperauthor{Martin Kalck}
\paperemail{m.kalck@ed.ac.uk}



\makepapertitle

\begin{abstract}
Leclerc recently studied certain Frobenius categories in connection with cluster algebra structures on coordinate rings of intersections of opposite Schubert cells. We show that these categories admit a description as Gorenstein projective modules over an Iwanaga-Gorenstein ring of virtual dimension at most two. This is based on a Morita type result for Frobenius categories. 
\end{abstract}

\section{Motivation}

\noindent
Let $G$ be a complex simple Lie group of type $Q=A, D$ or $E$ (eg $G=\mathrm{SL}_{n+1}(\mathbb{C})$ for $Q=A_n$) with Borel subgroup $B \subset G$ (eg $B=\{\text{upper triangular matrices}\}$) and Weyl group $W$ (eg $W \cong S_{n+1}$ given by permutation matrices).

For a Weyl group element $w \in W$ there are associated subvarieties $C_w$ (\emph{Schubert cell}) and $C^w$ (\emph{opposite Schubert cell}) in the flag variety $G/B$. On the other hand, there is a torsion pair $(\cc_w, \cc^w)$ in the category of finite dimensional modules over the preprojective algebra $\Pi:=\Pi(Q)$ and the categories $\cc_w$, $\cc^w$ are Frobenius and have projective generators (in fact, the latter statements may be deduced from Proposition \ref{P:Iyama}). These Frobenius categories were used by Gei{\ss}, Leclerc \& Schr\"oer to categorify cluster algebra structures on  coordinate rings of the corresponding (opposite) Schubert cells \cite{GLS}.

Let $v \in W$. The intersections $C_{v, w}:=C^v \cap C_w$ are known as \emph{open Richardson varieties} and have been studied by Kazhdan-Lusztig in connection with KL-polynomials. Generalizing the aforementioned work \cite{GLS}, Leclerc \cite{L} categorifies a cluster subalgebra of the coordinate rings of $C_{v, w}$ using the intersection $\cc_{v, w}$ of a torsion free part $\cc^v$ with a torsion part $\cc_w$ of two torsion pairs mentioned above. Under some finiteness assumptions he obtains a cluster algebra structure on the whole coordinate ring and he conjectures that this holds in general.

The subcategories $\cc_{v, w} \subseteq \mod \Pi$ inherit an exact structure which is again Frobenius.
\begin{aim}
Explain this in a more abstract setting and give equivalent descriptions of $\cc_{v, w}$.
\end{aim} 
\noindent
This is summarized in the following Proposition which is a special case of Proposition \ref{P:Iyama}.

\begin{proposition}\label{P:Main}
Let $\cc_{v, w}:=\cc_w \cap \cc^v \subseteq \mod \ \Pi$. Then 
\begin{itemize}
\item[(a)] $\cc_{v, w}$ is a Frobenius category with $\proj \cc_{v, w} = \add f_v t_w (\Pi) = \add t_w f_v (\Pi) =:\add P_{v, w}$. Where $t_u(-)$ denotes the torsion radical and $f_u(-):=(-)/t_u(-)$ for a torsion pair $(\cc_u, \cc^u)$.
\item[(b)]Ê$\cc_{v, w} \xrightarrow{\Hom_{\cc_{v, w}}(P_{v, w}, -)} \GP\bigl(\Pi_{v, w}\bigr)$ is an exact equivalence, where $\Pi_{v,w}:=\End_{\cc_{v, w}}(P_{v, w})$ is an Iwanaga-Gorenstein ring of virtual dimension at most two. 
\item[(c)] In particular, $\cc_{v, w}$ is equivalent to the subcategory of second syzygies of finite dimensional $\Pi_{v,w}$-modules.
\item[(d)] The functors $f_v$ and $t_w$ induce ring homomorphisms
$ \Pi_w:=\End_{\cc_w}(t_w(\Pi)) \to \Pi_{v, w}$ and 
$ \Pi^v:=\End_{\cc^v}(f_v(\Pi)) \to \Pi_{v, w}$. These 
are surjective if $\cc_v \subseteq \cc_w$. In turn, this condition is equivalent to $w=v'v$ with $l(w)=l(v') + l(v)$, called condition (P) in Leclerc \cite[5.1]{L}.
\item[(e)](see \cite[5.16]{BKT}) If condition (P) holds, then $\Pi_{v, w}$ is Morita equivalent to $\Pi_{v'}$.  Therefore, $\Pi_{v, w}$ has the same virtual dimension as $\Pi_{v'}$ which is at most $1$,  \cite{BIRS}.
\end{itemize}
\end{proposition}

\begin{remark}\label{R: BIRS}
Let $\Lambda_w:= \Pi/I_w$ be the algebra considered in \cite{BIRS}. Then there are algebra isomorphisms $ \Lambda_w \cong \Pi^{w_0 w^{-1}} \cong \Pi_{w^{-1}}^{\rm op}$, where $w_0$ denotes the longest Weyl group element.
\end{remark}

\section{A Morita type result for Frobenius categories} 

\begin{DefProp}
A two-sided Noetherian ring $R$ is called \emph{Iwanaga-Gorenstein}, if 
$
\injdim _R R < \infty$  and  $ \injdim R_R < \infty.$
It is well-known that this implies 
$
\injdim _R R =d= \injdim R_R 
$. We call $d=:\virdim R$ the \emph{virtual dimension} of $R$. 

In this case the category of \emph{Gorenstein-projective} $R$-modules
\begin{align*}
\GP(R):=\{ M \in \mod  R \mid \Ext^i_R(M, R)=0 \,\text{  for all }\, i >0 \}
\end{align*}
is a Frobenius category with subcategory of projective-injective objects $\proj R$. Equivalently, $\GP(R)$ is the subcategory of $d$-th syzygies of finitely generated $R$-modules
\begin{align*}
\GP(R) \cong \Omega^d(\mod R):=\{ \Omega^d(M)  \mid  M \in \mod R\}.  
\end{align*}
If $R$ is a local commutative Noetherian ring, Gorenstein projective $R$-modules are precisely maximal Cohen-Macaulay $R$-modules and $\injdim _R R=\krdim \, R$.
\end{DefProp}

\begin{aim} Characterize the categories of Gorenstein projective modules $\GP(R)$ over Iwanaga-Gorenstein  rings $R$ among all Frobenius categories.
\end{aim}

\begin{notation}
For an additive category $\cb$, we denote by $\mod \cb$ the category of finitely presented contravariant additive functors $\cb \to \mathrm{Ab}$.
\end{notation}

We first list properties of the categories $\ce:=\GP(R)$ for $R$ Iwanaga-Gorenstein.
\begin{itemize}
\item[(i)] $\proj \ce=\add P \, (=\proj R)$ for some $P \in \ce$ and $\End_\ce(P)\, (\cong \End_R(R))$ is two-sided noetherian.
\item[(ii)] $\ce$ is idempotent complete (since $\ce \subseteq \mod R$ closed under direct summands).
\item[(iii)] $\ce$ is Frobenius (use exact duality $\Hom_R(-, R) \colon \GP(R) \to \GP(R^{\rm op})$).
\item[(iv)] $\ce$ has weak kernels and cokernels (use Auslander-Buchweitz approximation).
\item[(v)] $\gldim \mod  \ce$, $\gldim \mod  \ce^{\rm op} \leq n \, (= \max\{2, \injdim R \}). $ 
\end{itemize}

The following result may be interpreted as an analogue of Morita theory for Frobenius categories. The implication (b) $\Rightarrow$ (a) is well-known. The converse is the special case $\proj \ce=\add \cp, \cm=\ce$ of {\cite[2.8]{KIWY}}, which is due to Iyama and inspired by a stable version of Dong Yang and the author \cite[2.15]{KIWY}. 

\begin{proposition} Let $\ce$ be an exact category and let $P \in \ce$. TFAE
\begin{itemize}
\item[(a)] $\ce$ and $P$ satisfy the conditions (i)-(v) above.
\item[(b)] Set $R=\End_\ce(P)$. $\Hom_\ce(P, -)\colon \ce \to \GP(R)$ is an exact equivalence and $R$ is Iwanaga-Gorenstein with  $\virdim R \leq \gldim \mod  \ce$.
\end{itemize}
\end{proposition}

\section{From pairs of torsion pairs to Frobenius categories}

\begin{notation}
Let  $(\ct, \cf)$ be a torsion pair in an abelian category $\ca$. In particular, there is a short exact sequence $0 \to t(X) \to X \to f(X) \to 0$ for all $X$ in $\ca$. This gives rise to functors $t\colon \ca \to \ct$ and $f\colon \ca \to \cf$, which are right (respectively left) adjoint to the canonical inclusions. 
\end{notation}

\begin{proposition} \label{P:Iyama}
Let $\ca$ be an abelian category with torsion pairs $(\ct_1, \cf_1)$ and $(\ct_2, \cf_2)$ and set $\cc_{12}:=\ct_1 \cap \cf_2$. Then the following statements hold:
\begin{itemize}
\item[(a)] $\cc_{12}$ is extension closed and idempotent complete, since $\ct_1$ and $\cf_2$ are. In particular, $\cc_{12}$ inherits a natural exact structure from $\ca$.
\item[(b)] $\cc_{12}$ has kernels and cokernels. In other words, $\cc_{12}$ is a preabelian category. In particular, the categories of finitely presented additive functors $\mod \, \cc_{12}$ and $\mod \, \cc_{12}^{\rm op}$ are abelian and have global dimension at most $2$.

For example, the composition of the canonical inclusions
\[
t_1(\ker f) \hookrightarrow \ker f \hookrightarrow X \xrightarrow{f} Y
\]
is a kernel of $f$. Here $\ker f$ denotes the kernel of $f$ in $\ca$.
\item[(c)] If $\ct_1$ has enough projectives and $\cf_2$ has enough injectives, then $\cc_{12}$ has enough injectives ($=t_1(\mathrm{inj} \ \cf_2)$) and projectives ($=f_2(\proj \ct_1)$).
\item[(d)] If additionally $\Ext^1_{\cc_{12}}(X, Y)=0 \Leftrightarrow \Ext^1_{\cc_{12}}(Y, X)=0$, then $\cc_{12}$ is Frobenius. 
For example, this is satisfied if $\ul{\ca}$ or $\cd^b(\ca)$ are $2$-Calabi-Yau. This in turn is known to hold for $\ca=\mathrm{fdmod}(\widehat{\Pi(Q)})$, where $Q$ is a quiver without loops and $\widehat{\Pi(Q)}$ is the $\mathrm{m}$-adic completion of its preprojective algebra, where $\mathrm{m}$ denotes the ideal generated by all arrows.
\item[(e)] Assume additionally that $\proj \ \ct_1=\add P$ and $\mathrm{inj} \ \cf_2=\add I$, then $\proj \cc_{12}=\add f_2(P) = \add t_1(I)$. We assume that $\Pi_{12}:=\End_{\cc_{12}}(f_2(P))$ is two-sided noetherian. Then there is an exact equivalence
\begin{align*}
\cc_{12} \xrightarrow{\Hom_{\cc_{12}}(f_2(P), -)} \GP\bigl(\Pi_{12}\bigr),
\end{align*}
and $\Pi_{12}$ is Iwanaga-Gorenstein of virtual dimension at most $2$.
\item[(f)] In the situation of (e) the functors $f_2$ and $t_1$ induce ring homomorphisms $\varphi_{2}\colon \End_{\ct_1}(P) \to \Pi_{12}$ and $\tau_1\colon \End_{\cf_2}(I) \to \Pi_{12}$ with kernels given by the ideals of morphisms factoring over $t_2(P)$ and $f_1(I)$, respectively. The ring homomorphisms are surjective if $\ct_2 \subseteq \ct_1$. In Example \ref{E:Leclerc}, $\varphi_2$ is injective but not surjective.

\end{itemize}

\end{proposition}

\begin{remark}
 This is an analogue of Buan, Iyama, Reiten \& Scott's \cite{BIRS} dual description of Gei{\ss}, Leclerc \& Schr{\"o}er's categories $\cc_w$ \cite{GLS} as categories of submodules of projective modules over the algebra $\Lambda_w$, see also \cite[Theorem 2.8]{GLS}. Since $\Lambda_w$ is Iwanaga-Gorenstein of virtual dimension $1$, Gorenstein projective modules are first syzygies, which in turn are just submodules of projective modules. See also \cite[Section 6]{KIWY} for a further discussion.\end{remark}

\section{Examples, remarks and questions}



\begin{example}\label{E:Leclerc}
We consider the situation of \cite[3.16]{L}, i.e. $Q$ is of type $A_3$, $w=s_1s_3s_2s_1s_3$ and $v=s_2$. Then $\varphi_2 \colon \Pi_w:=\End_{\cc_w}(t_w(\Pi)) \to \Pi_{v, w}$ is injective and its cokernel in the category of vectorspaces is isomorphic to $\mathbb{C}$. Moreover, $\Pi_{v, w}$ is the Auslander algebra of the preprojective algebra of type $A_2$ and therefore is of global (and virtual) dimension $2$.
\end{example}

\begin{remark}[Duality] Let $Q$ be a Dynkin quiver and let $D:=\Hom_k(-, k)$ be the standard duality. It is well-known that there is an algebra isomorphism $\psi \colon \Pi \cong \Pi^{\rm op}$, which gives rise to a duality $\Phi \colon  \mod \Pi \xrightarrow{D}  \mod \Pi^{\rm op} \xrightarrow{\psi_*}  \mod \Pi$. Using the notation in Leclerc \cite[\textsection 3.2]{L}, one can check that $\Phi(P_{v, w}) \cong P_{ w_0^{-1}w, w_0 v}$ holds, where $w_0$ denotes the longest Weyl group element. In particular, $\Phi$ induces an algebra isomorphism $\Pi_{v, w} \cong \Pi_{ w_0^{-1}w, w_0 v}^{\rm op}$.  Thus $\Pi^v \cong \Pi_{v, w_0} \cong \Pi_{\mathrm{id}, w_0v}^{\rm op} \cong \Pi_{w_0v}^{\rm op}$ for the algebras appearing in Proposition \ref{P:Main} (d).
\end{remark}

\begin{open}
Give a 'combinatorial description' of $\Pi_{v, w}$, eg as quiver with relations.
\end{open}

\begin{remark}
 The number of isoclasses of indecomposable projective $\Pi_{v, w}$-modules seems to be unknown in general. It is not always bounded above by $\mid \! \! Q_0 \! \! \mid$, see Example \ref{E:Leclerc}.
\end{remark}

\begin{question}[Leclerc]
How does the virtual dimension of $\Pi_{v, w}$ depend on $Q, v, w$ and  
(how) is this number related to the geometry of the open Richardson variety $C_{v, w}$?
\end{question}

\begin{partial answer} By Remark \ref{R: BIRS} and \cite{BIRS}, $\virdim \Pi^v, \Pi_w \leq 1$. They are zero iff $\cc^v$ (respectively, $\cc_w$) are exact abelian subcategories of $\mod \Pi$, which are then equivalent to $\mod \Pi/e$ ($e \in \Pi$ idempotent). Thus if $\virdim \Pi^v, \Pi_w =0$, then $\virdim \Pi_{v, w}=0$  (since $\cc_{v, w}$ is abelian). If one of $\Pi^v$ and $\Pi_w$ has virtual dimension zero, then $\cc_{v,w}$ is the torsion (or torsion-free) part of a torsion pair in $\mod \Pi/e$. By Mizuno \cite{M} and \cite{BIRS}, $\Pi_{v,w} \cong (\Pi/e)^{v'}$ or $ \cong (\Pi/e)_{w'}$ and is therefore of virtual dimension $\leq 1$. Also $\gldim \Pi_v=n \leq 1$ (or $\gldim \Pi^w =m \leq 1$) implies $\gldim \Pi_{v,w} \leq \min\{n, m\}$. If both $\Pi_v, \Pi^w$ have infinite global dimension and virtual dimension $1$, then virtual dimensions $0, 1, 2$ occur for $\Pi_{v, w}$.
\end{partial answer}

\begin{remark}[Commutativity]
It follows from work of Mizuno \cite{M}, that all torsion pairs in $\mod \Pi$ are of the form $(\cc_w, \cc^w)$ for some Weyl group element $w$. In particular, there are only finitely many torsion pairs, which is very surprising given the size of $\mod \Pi$. The explicit description of the associated functors $t_w$ and $f_w$ (see eg Leclerc \cite[\textsection 3.2]{L}) shows that $f_v t_w(M) \cong  t_w f_v(M)$ for Weyl group elements $v, w \in W$ and $M \in \mod \Pi$. This seems very unusual for a pair of torsion pairs in general abelian categories and fails already for $\mod U_2(k)$, where $U_2(k)$ denotes the ring of $2 \times 2$ upper triangular matrices.
\end{remark}

\noindent
\emph{Acknowledgement}. This material grew out of a discussion with Bernard Leclerc and Henning Krause after a talk of Leclerc on the results of \cite{L} in spring 2014. I am grateful to Bernard Leclerc for his inspiring work, his interest and discussions. Moreover, I would like to thank Henning Krause for insisting that these considerations might be of interest and Osamu Iyama for very inspiring discussions. In particular, I learned parts of Proposition \ref{P:Iyama} from him. I am very grateful to Michael Wemyss for lots of stimulating questions on this topic. I also had fruitful discussions with Sergio Estrada, Mikhail Gorsky, Frederik Marks, Yuya Mizuno \& Milen Yakimov. 
Thanks to the organizers of this conference for the opportunity to present this work and to the participants for their interest and questions. I am grateful to EPSRC for financial support (EP/L017962/1).

\end{document}